\documentclass[11pt]{article}
\usepackage{amsmath}
\usepackage{amssymb}
\usepackage{amsthm}
\usepackage{url}
\newcommand\nonu{\nonumber}
\newcommand\sLP{\\[\smallskipamount]}
\newcommand\sPP{\\[\smallskipamount]\indent}
\newcommand\bPP{\\[\bigskipamount]\indent}
\newcommand\CC{\mathbb{C}}
\newcommand\RR{\mathbb{R}}
\newcommand\ZZ{\mathbb{Z}}
\newcommand\al\alpha
\newcommand\be\beta
\newcommand\ga\gamma
\newcommand\de\delta
\newcommand\la\lambda
\newcommand\tha\theta
\newcommand\Ga{\Gamma}
\newcommand\De{\Delta}
\newcommand\La{\Lambda}
\newcommand\Om{\Omega}
\newcommand\half{\frac12}
\newcommand\thalf{\tfrac12}
\newcommand\iy\infty
\newcommand\wt{\widetilde}
\newcommand\LHS{left-hand side}
\newcommand\RHS{right-hand side}
\renewcommand\Re{\operatorname{Re}}
\newcommand\const{\mathrm{const.}\,}
\newcommand{\hyp}[5]{\,\mbox{}_{#1}F_{#2}\!\left(
  \genfrac{}{}{0pt}{}{#3}{#4};#5\right)}
\newcommand{\qhyp}[5]{\,\mbox{}_{#1}\phi_{#2}\!\left(
  \genfrac{}{}{0pt}{}{#3}{#4};#5\right)}

\numberwithin{equation}{section}
\newtheorem{theorem}{Theorem}[section]
\newtheorem{proposition}[theorem]{Proposition}

\newtheorem{Definition}[theorem]{Definition}

\newtheorem{Remark}[theorem]{Remark}
\newenvironment{remark}{\begin{Remark}\rm}{\end{Remark}}
\newtheorem{Example}[theorem]{Example}
\newenvironment{example}{\begin{Example}\rm}{\end{Example}}
\newcommand\Proof{\noindent{\bf Proof}\quad}

\begin{document}
\title{Quadratic transformations for orthogonal polynomials\\
in one and two variables}
\author{Tom H. Koornwinder}
\date{Dedicated to Masatoshi Noumi on the occasion of his sixtieth
birthday}
\maketitle
\begin{abstract}
We discuss quadratic transformations for orthogonal polynomials
in one and two variables. In the one-variable case
we list many (or all) quadratic transformations between families
in the Askey scheme or $q$-Askey scheme. In the two-variable case
we focus, after some generalities, on the polynomials associated
with root system $BC_2$, i.e., $BC_2$-type Jacobi polynomials
if $q=1$ and Koornwinder polynomials in two variables in the $q$-case.

\end{abstract}
% until 74
\section{Introduction}
Whenever we have a system of orthogonal polynomials $\{p_n\}$
in one variable with respect to an even orthogonality measure $\mu$
on $\RR$,
then we can write $p_{2n}(x)=q_n(x^2)$, $p_{2n+1}(x)=x\,r_n(x^2)$
with $\{q_n\}$ and $\{r_n\}$ systems of orthogonal polynomials on
$\RR_{\ge0}$ with respect to orthogonality measures which are
immediately obtained from $\mu$. These mappings from
$\{p_n\}$ to $\{q_n\}$ and $\{r_n\}$ are called
quadratic transformations.
For quite some multi-parameter families of orthogonal polynomials
in the Askey scheme and the $q$-Askey scheme such quadratic
transformations can be given explicitly. Very well-known
are the quadratic transformations for Jacobi polynomials
connecting $\big\{P_n^{(\al,\al)}\big\}$ with
$\big\{P_n^{(\al,\pm\half)}\big\}$.
Since all such
polynomials can be expressed as ($q$-)hypergeomnetric functions,
their quadratic transformations are equivalent to certain quadratic
transformations for terminating \mbox{($q$-)hypergeometric} functions.

The first aim of this paper, in Section 2,
is to survey many (maybe all) instances of
quadratic transformations in the ($q$-)Askey scheme, and how they are
related by the limit arrows in those schemes.
While the quadratic transformations for Askey-Wilson polynomials
were already given in the Memoir \cite{1} by Askey \& Wilson,
some of the other quadratic transformations given below may
occur here for the
first time, in particular the ones on the discrete side of the
($q$-)Askey scheme.

Quadratic transformations occur also for orthogonal polynomials
in several variables as soon as the orthogonality measure is
invariant under the transformation $x_1\mapsto-x_1$ of the first
variable~$x_1$. This sounds like a trivial generalization
of the one-variable case,
but this reflection map already takes some unexpected form
when we look for quadratic transformations within
multi-parameter families of special orthogonal polynomials
in two variables. For the systems associated with root system $BC_2$
the deeper explanation for the existence of the quadratic
transformations is the isomorphism between the root systems $B_2$ and
$C_2$, both of which are contained in $BC_2$.

These quadratic transformations in the two-variable case
will be discussed in Section 3.
For $BC_2$-type Jacobi polynomials they go back to
Sprinkhuizen-Kuyper \cite{4}, while they may be new for
Koornwinder polynomials. We will also argue that quadratic
transformations for orthogonal polynomials associated with $BC_n$
cannot occur if $n>2$, at least not in the simple form as for
$n=1$ and 2.

The paper concludes in Section 4 with a discussion how quadratic
transformations can be helpful as heuristics for extending results
to a larger realm of parameters, and with mentioning some
possible work which would be a natural follow-up of this paper.
\paragraph{Conventions}
For definition and notation of
hypergeometric and $q$-hypergeometric series
see~\cite{3}. Throughout we will assume
that $0<q<1$.
\section{The one-variable case}
\subsection{Ordinary polynomials}
Let $\{p_n(x)\}$ be a system of monic orthogonal polynomials on $\RR$
which are orthogonal with respect to an even (nonnegative) weight function
$w(x)=w(-x)=v(x^2)$. Then $p_n(-x)=(-1)^n p_n(x)$.
Put
\begin{equation}
q_n(x^2):=p_{2n}(x),\quad
r_n(x^2):=x^{-1} p_{2n+1}(x).
\label{1}
\end{equation}
Then (see \cite[Ch.~1, \S8]{10})
$\{q_n(x)\}$ and $\{r_n(x)\}$ are systems of monic orthogonal
polynomials on $[0,\iy)$:\\
$\bullet$\; the $q_n$ with respect to weight function $x^{-\half} v(x)$,\\
$\bullet$\; the $r_n$ with respect to weight function $x^\half v(x)$.
\sPP
Note that from \eqref{1} we have, for any $x_0\in\RR$ on which the $p_n$
do not vanish, that
\begin{equation}
\frac{q_n(x^2)}{q_n(x_0^2)}=\frac{p_{2n}(x)}{p_{2n}(x_0)}\,,\quad
\frac{r_n(x^2)}{r_n(x_0^2)}=\frac{x_0\,p_{2n+1}(x)}{x\,p_{2n+1}(x_0)}\,.
\label{2}
\end{equation}
The identities \eqref{2} remain valid for arbitrary normalizations of the
$p_n$, $q_n$, $r_n$.

As a slight variant of the above,
let $\{p_n(x)\}$ be a system of orthogonal polynomials on $[-1,1]$
which are orthogonal with respect to an even weight function
$w(x)=w(-x)=v(2x^2-1)$. Let $x_0\in\RR$ such that $p_n(x_0)\ne0$ for all $n$.
Let $q_n(x)$ and $r_n(x)$ be polynomials of degree $n$ such that
\begin{equation*}
\frac{q_n(2x^2-1)}{q_n(2x_0^2-1)}=\frac{p_{2n}(x)}{p_{2n}(x_0)},\quad
\frac{r_n(2x^2-1)}{r_n(2x_0^2-1)}=
\frac{x_0\,p_{2n+1}(x)}{x\,p_{2n+1}(x_0)}\,.
\end{equation*}
Then $\{q_n(x)\}$ and $\{r_n(x)\}$ are systems of orthogonal
polynomials on $[-1,1]$:\\
$\bullet$\; the $q_n$ with respect to weight function
$(1+x)^{-\half} v(x)$,\\
$\bullet$\; the $r_n$ with respect to weight function $(1+x)^\half v(x)$.
\begin{example}
\label{4}
Jacobi polynomials
\begin{equation*}
P_n^{(\al,\be)}(x)=(-1)^nP_n^{(\be,\al)}(-x):=
\frac{(\al+1)_n}{n!}\,\hyp21{-n,n+\al+\be+1}{\al+1}{\thalf(1-x)}
\end{equation*}
are orthogonal on $[-1,1]$ with weight funnction $(1-x)^\al(1+x)^\be$\quad
($\al,\be>-1$).
So we have quadratic transformations (see \cite[Theorem 4.1]{17})
\begin{equation}
\frac{P_{2n}^{(\al,\al)}(x)}{P_{2n}^{(\al,\al)}(1)}
=\frac{P_n^{(\al,-\half)}(2x^2-1)}{P_n^{(\al,-\half)}(1)}\,,\quad
\frac{P_{2n+1}^{(\al,\al)}(x)}{P_{2n+1}^{(\al,\al)}(1)}
=\frac{x P_n^{(\al,\half)}(2x^2-1)}{P_n^{(\al,\half)}(1)}\,.
\label{10}
\end{equation}
\end{example}

\begin{example}
Laguerre polynomials
\begin{equation*}
L_n^\al(x):=\frac{(\al+1)_n}{n!}\,\hyp11{-n}{\al+1}x
\end{equation*}
are orthogonal on $[0,\iy)$ with
weight function $x^\al e^{-x}$\quad($\al>-1$), while
Hermite polynomials 
\begin{equation*}
H_n(x):=(2x)^n\,\hyp20{-\thalf n,-\thalf(n-1)}-{-x^{-2}}
\end{equation*}
are orthogonal on $(-\iy,\iy)$
with weight function $e^{-x^2}$. So we have quadratic trnansformations
(see \cite[(5.6.1)]{17})
\begin{equation}
\frac{H_{2n}(x)}{H_{2n}(0)}=\frac{L_n^{-\half}(x^2)}{L_n^{-\half}(0)}\,,
\qquad
\frac{H_{2n+1}(x)}{H_{2n+1}'(0)}
=\frac{x L_n^{\half}(x^2)}{L_n^{\half}(0)}\,.
\label{53}
\end{equation}
These are limit cases of \eqref{10} by the limits
\cite[(9.8.16), (9.8.18)]{13}.
\end{example}

\begin{remark}
\label{41}
In connection with \eqref{1} and \eqref{2} we had weight functions
$w(x)=w(-x)=v(x^2)$. Then $d\mu(x):=w(x)\,dx$ is an even measure on $\RR$
and $d\nu(x):=2x^{-\half}v(x)\,dx$ is the pushforward measure
$\nu=\phi_*\mu$ on $\RR_{\ge0}$ with
$\phi\colon x\mapsto x^2\colon\RR\to\RR_{\ge0}$.
In general, the quadratic transformations \eqref{1}, \eqref{2} remain true
if the $p_n$ are orthogonal polynomials with respect
to a (positive) even measure on $\RR$, the $q_n$ are orthogonal
with respect to the measure $\nu=\phi_*\mu$ on $\RR_{\ge0}$, i.e.,
\[
\int_{\RR_{\ge0}} p(y)\,d\nu(y)=\int_\RR p(x^2)\,d\mu(x)
\quad\mbox{for all polynomials $p$},
\]
and the $r_n$ are orthogonal with respect to the measure $x\,d\nu(x)$
on $\RR_{\ge0}$. Similar remarks will apply to other quadratic transformations.
This becomes in particular relevant in examples involving discrete mass
points or $q$-integrals.
\end{remark}

\subsection{Symmetric Laurent polynomials}
As a further variant of the above, with $w(x)$ a weight function on $[-1,1]$,
we substitute
$x=\thalf(z+z^{-1})$ so that $z$ runs from $-1$ to 1 on the upper half unit
circle if $x$ runs from $-1$ to 1 on the interval $[-1,1]$.
Let $\De(z)$ be a real-valued
weight function on the upper half unit circle such that
\begin{equation*}
w(x)=w\big(\thalf(z+z^{-1})\big)=\frac{2i\De(z)}{z-z^{-1}}\,.
\end{equation*}
Then
\begin{equation*}
\int_{-1}^1 f(x)\,w(x)\,dx=
i^{-1}\int_C f\big(\thalf(z+z^{-1})\big)\,\De(z)\,\frac{dz}z\,,
\end{equation*}
where the contour $C$ is the upper half unit circle starting at 1 and
ending at $-1$.
Now suppose that $\De(z)=\De(-z^{-1})$ and put
\begin{equation}
\wt\De(z^2):=\De(z)=\De(-z^{-1})\qquad(|z|=1,\;0\le\arg z\le\pi/2).
\label{3}
\end{equation}
Equivalently, $w(x)=w(-x)$. As before, put
\begin{equation*}
v(2x^2-1):=w(x)=w(-x).
\end{equation*}
Then
\begin{equation*}
w\big(\thalf(z+z^{-1})\big)
=v\big(\thalf(z^2+z^{-2})\big)
=\frac{2i}{z-z^{-1}}\,\wt\De(z^2).
\end{equation*}
Hence
\begin{align*}
\big(1+\thalf(z^2+z^{-2})\big)^{-\half}v\big(\thalf(z^2+z^{-2})\big)
&=2^\half\,\frac{2i}{z^2-z^{-2}}\,\wt\De(z^2),\\
\big(1+\thalf(z^2+z^{-2})\big)^{\half}v\big(\thalf(z^2+z^{-2})\big)
&=2^{-\half}(1+z^2)(1+z^{-2})\,\frac{2i}{z^2-z^{-2}}\,\wt\De(z^2).
\end{align*}
Thus, with $x=\thalf(z+z^{-1})$,
\begin{align*}
(1+x)^{-\half}v(x)&=2^\half\,\frac{2i}{z-z^{-1}}\,\wt\De(z),\\
(1+x)^\half v(x)&=2^{-\half}(1+z)(1+z^{-1})\,\frac{2i}{z-z^{-1}}\,\wt\De(z).
\end{align*}
We arrive at the following result.
Let $\{\wt p_n(z)\}$ be a system of symmetric (i.e., invariant under
$z\to z^{-1}$) Laurent polynomials
which are orthogonal on $C$ with respect to the measure
$\De(z) z^{-1}dz$,
where $\De$ satisfies \eqref{3}.
Let $z_0\in\CC$ such that $p_n(z_0)\ne0$ for all $n$.
Let $\wt q_n(z)$ and $\wt r_n(z)$ be symmetric Laurent polynomials of
degree $n$ such that
\begin{equation*}
\frac{\wt q_n(z^2)}{\wt q_n(z_0^2)}=\frac{\wt p_{2n}(z)}{\wt p_{2n}(z_0)},\quad
\frac{\wt r_n(z^2)}{\wt r_n(z_0^2)}=
\frac{(z_0+z_0^{-1})\wt p_{2n+1}(z)}{(z+z^{-1})\wt p_{2n+1}(z_0)}\,.
\end{equation*}
Then $\{\wt q_n(z)\}$ and $\{\wt r_n(z)\}$ are systems of
symmetric orthogonal Laurent
polynomials on $C$:\\
$\bullet$\; the $\wt q_n$ with orthogonality measure
$\wt\De(z) z^{-1}dz$,\\
$\bullet$\; the $\wt r_n$ with orthogonality measure
$(1+z)(1+z^{-1}) \wt\De(z) z^{-1}dz$.

If we go back to Example \ref{4} then, with the above notation and up
to constant factors,
\begin{align*}
&\wt p_n(z)=P_n^{(\al,\al)}\big(\thalf(z+z^{-1})\big),\quad
\wt q_n(z)=P_n^{(\al,-\half)}\big(\thalf(z+z^{-1})\big),\quad
\wt r_n(z)=P_n^{(\al,\half)}\big(\thalf(z+z^{-1})\big),\\
&\De(z)=(2-z^2-z^{-2})^{\al+\half},\quad
\wt\De(z)=(2-z-z^{-1})^{\al+\half}.
\end{align*}
\begin{example}
Recall Askey-Wilson polynomials \cite{1}, \cite[\S14.1]{13}, which we write as monic
symmetric Laurent polynomials:
\begin{multline}
P_n(z)=P_n(z;a,b,c,d\,|\,q):=\frac1{(abcdq^{n-1};q)_n}\,
p_n\big(\thalf(z+z^{-1});a,b,c,d\,|\,q\big)\\
=\frac{(ab,ac,ad;q)_n}{a^n (abcdq^{n-1};q)_n}\,
\qhyp43{q^{-n},q^{n-1}abcd,az,az^{-1}}{ab,ac,ad}{q,q}.
\label{48}
\end{multline}
Here $P_n(z)$ is invariant under permutations of the parameters
$a,b,c,d$. Observe that
\begin{equation*}
P_n(a;a,b,c,d\,|\,q)=\frac{(ab,ac,ad;q)_n}{a^n (abcdq^{n-1};q)_n}\,,
\quad
p_n(\thalf(a+a^{-1});a,b,c,d\,|\,q)=\frac1{(abcdq^{n-1};q)_n}\,,
\end{equation*}
\begin{equation}
\frac{p_n(\thalf(z+z^{-1}))}{p_n(\thalf(a+a^{-1}))}=
\qhyp43{q^{-n},q^{n-1}abcd,az,az^{-1}}{ab,ac,ad}{q,q}.
\label{33}
\end{equation}

Assume that $a,b,c,d$ have absolute value $\le1$ but
do not have pairwise products equal to~1, and that non-real parameters
occur in complex conjugate pairs.
The polynomials $P_n(z)$ are orthogonal
on the upper half unit circle $C$ with
respect to the orthogonality measure
$\De(z)\,z^{-1}\,dz$, where
\[
\De(z)=\De_+(z)\De_+(z^{-1}),\qquad
\De_+(z)=\De_+(z;a,b,c,d\,|\,q):=\frac{(z^2;q)_\iy}
{(az,bz,cz,dz\,|\,q)_\iy}\,.
\]
Since
\begin{equation*}
\De(z;a,b,-a,-b\,|\,q)=\De(z^2;a^2,b^2,-1,-q\,|\,q^2)
=\frac{\De(z;a^2,b^2,-q,-q^2\,|\,q^2)}{(1+z^2)(1+z^{-2})}\,,
\end{equation*}
we have:
\begin{align}
P_{2n}(z;a,b,-a,-b\,|\,q)&=P_n(z^2;a^2,b^2,-1,-q\,|\,q^2),
\label{5}\\
P_{2n+1}(z;a,b,-a,-b;q)&=(z+z^{-1}) P_n(z^2;a^2,b^2,-q,-q^2\,|\,q^2),
\label{6}
\end{align}
or, in the normalization \eqref{33},
\begin{align}
\frac{p_{2n}(x;a,b,-a,-b\,|\,q)}
{p_{2n}(\thalf(a+a^{-1});a,b,-a,-b\,|\,q)}
&=\frac{p_n(2x^2-1;a^2,b^2,-1,-q\,|\,q^2)}
{p_n(\thalf(a^2+a^{-2});a^2,b^2,-1,-q\,|\,q^2)}\,,
\label{34}\sLP
\frac{p_{2n+1}(x;a,b,-a,-b;q)}{p_{2n+1}(\thalf(a+a^{-1});a,b,-a,-b;q)}
&=\frac{2x\,p_n(2x^2-1;a^2,b^2,-q,-q^2\,|\,q^2)}
{(a+a^{-1})\,p_n(\thalf(a^2+a^{-2});a^2,b^2,-q,-q^2\,|\,q^2)}\,.
\label{35}
\end{align}

Formula \eqref{34} is given by Askey \& Wilson \cite[Section 3.1]{1}
in terms of $q$-hypergeometric functions, but similarly derived as above
($a,b$ below different from $a,b$ above):
\begin{equation}
\qhyp43{a^2,qb^2,c,-d}{qab,-qab,cd}{q,q}
=\qhyp43{a^2,qb^2,c^2,d^2}{q^2a^2b^2,cd,qcd}{q^2,q^2}
\label{7}
\end{equation}
when both sides terminate.
The identity \eqref{7} can also be obtained from
Singh \cite[(22)]{2} (see also \cite[(3.10.11)]{3}) by applying
Sears' transformation \cite[(2.10.4)]{3}.

While we arrived at \eqref{7} in the terminating case $a=q^{-n}$,
the identity
holds also in the terminating case $c=q^{-n}$. Then a resulting
identity for Askey-Wilson polynomials is
\begin{equation}
P_n(z;a,b,q^\half,-q^\half\,|\,q)
=P_n(z;a,qa,b,qb\,|\,q^2).
\label{8}
\end{equation}
This relates two different ways of writing continuous $q$-Jacobi polynomials
as Askey-Wilson polynomials, see \cite[(4.20)]{1} or \cite[(7.5.26)]{3}.
Formula \eqref{8} also follows by observing that
\begin{equation*}
\De(z;a,b,q^\half,-q^\half\,|\,q)=\De(z;a,qa,b,qb\,|\,q^2).
\end{equation*}

The quadratic transformation \eqref{35} can be written in terms of
$q$-hypergeometric functions as
\begin{equation}
\qhyp43{a^2,qb^2,c,-d}{qab,-qab,cd}{q,q}
=\frac{c-d}{1-cd}
\qhyp43{qa^2,q^2b^2,c^2,d^2}{q^2a^2b^2,q^2cd,qcd}{q^2,q^2}
\label{9}
\end{equation}
when both series terminate. For $c=q^{-n}$ the resulting identity
for Askey-Wilson polynomials is again \eqref{8},
with $a$ and $qa$ interchanged in the parameter list on the \RHS.
Hence, if $c=q^{-n}$ then \eqref{9}
follows from \eqref{7} by applying  Sears' transformation to the \RHS\
of \eqref{7}. 

With $ab=q^{\al+1}$ formulas \eqref{34}, \eqref{35}
give a two-parameter $q$-analogue of \eqref{10}.
Indeed if $a=a_q$, $b=b_q$ in \eqref{34}, \eqref{35} such that
$a_qb_q=q^{\al+1}$ and $a_q\to 1$ as $q\uparrow1$ then the quadratic
transformations \eqref{34}, \eqref{35} have the quadratic transformations
\eqref{10} as limits for $q\uparrow1$.

The quadratic transformations \eqref{34}, \eqref{35} remain valid
for less constrained parameter values by analytic continuation.
In the case of orthogonality involving additionally a finite number
of mass points (see \cite[(14.1.3)]{13}) we may still give a proof
of \eqref{34}, \eqref{35} by orthogonality in view of Remark \ref{41}.

There are various noteworthy special cases of the quadratic transformations
\eqref{34}, \eqref{35}. For $b=q^\half a$ we get continuous $q$-Jacobi
polynomials on the \LHS s and continuous $q$-ultraspherical polynomials
on the \RHS s.
For $b=0$ we get Al-Salam-Chihara
polynomials on the \LHS s and continuous dual $q$-Hahn polynomials
on the \RHS s.
For $a=b=0$ we get continuous $q$-Hermite polynomials 
on the \LHS s and Al-Salam-Chihara polynomials (in this context also
called continuous $q$-Laguerre polynomials) on the \RHS s.
See \cite[Ch.~14]{13} for details about the mentioned families of
orthogonal polynomials.
\end{example}

\subsection{Further examples of quadratic transformations in the
$q$-Askey scheme}
First we discuss some limit cases of the quadratic transformations
\eqref{34}, \eqref{35} for Askey-Wilson polynomials, where we stay
in the continuous part of the $q$-Askey scheme.
\begin{example}
For big $q$-Jacobi polynomials \cite[\S14.5]{13}
\begin{equation*}
P_n(x;a,b,c,d;q)=P_n(qac^{-1}x;a,b,-ac^{-1}d;q)
:=\qhyp32{q^{-n},q^{n+1}ab,qac^{-1}x}{qa,-qac^{-1}d}{q,q}.
\end{equation*}
and little $q$-Jacobi polynomials \cite[\S14.12]{13}, \cite[\S2.4]{12}
\begin{equation}
p_n(x;a,b;q):=\qhyp21{q^{-n},q^{n+1}ab}{qa}{q,qx},\quad
p_n(q^{-1}b^{-1};a,b;q)=
\frac{(-1)^n(qb;q)_n}{q^{\half n(n+1)}b^n(qa;q)_n}
\label{36}
\end{equation}
there are the quadratic transformations
\begin{align}
P_{2n}(x;a,a,1,1;q)&=\frac{p_n(x^2;q^{-1},a^2;q^2)}{p_n((qa)^{-2};q^{-1},a^2;q^2)}\,,
\label{37}\\
P_{2n+1}(x;a,a,1,1;q)&=\frac{qa\,x\,p_n(x^2;q,a^2;q^2)}{p_n((qa)^{-2};q,a^2;q^2)}\,.
\label{38}
\end{align}
These were earlier given in \cite[(2.48), (2.49)]{12}.
They are limit cases of \eqref{34}, \eqref{35}
by the limit formulas \cite[(6.2), (6.4)]{11}.

The orthogonality relations for big and little $q$-Jacobi polynomials
are given by $q$-integrals. In view of Remark \ref{41}
the quadratic transformations \eqref{37}
and \eqref{38} can be obtained in a straightforward way by comparing
the $q$-weights for the polynomials involved. The relevant observation is
that, with $w(x)=v(x^2)$ and polynomials $p$, we have
\begin{multline*}
\int_0^1 p(x)\,x^{-\half} v(x)\,d_{q^2}x
=(1-q^2)\sum_{k=0}^\iy p(q^{2k})\,v(q^{2k})\,q^k\\
=(1-q^2)\sum_{k=0}^\iy p((q^k)^2)\,w(q^k)\,q^k
=(1+q)\int_0^1 p(x^2)\,w(x)\,d_qx.
\end{multline*}
\end{example}

\begin{example}
For discrete $q$-Hermite I polynomials \cite[\S14.28]{13}
\begin{equation*}
h_n(x;q):=q^{\half n(n-1)} \qhyp21{q^{-n},x^{-1}}0{q,-qx}
\end{equation*}
and the little $q$-Laguerre polynomials (or Wall polynomials)
$p_n(x;a;q)=p_n(x;a,0;q)$ \cite[\S14.20]{13} (\eqref{36}
with $b=0$) there are the quadratic transformations
\begin{align}
h_{2n}(x;q)
&=(-1)^n q^{n(n-1)}(q;q^2)_n,p_n(x^2;q^{-1};q^2),
\label{39}\\
h_{2n+1}(x;q)
&=(-1)^n q^{n(n-1)}(q^3;q^2)_n\,x\,p_n(x^2;q;q^2).
\label{40}
\end{align}
These are limit cases of \eqref{37}, \eqref{38}
by the limit formula \cite[\S14.5]{14}
\begin{equation*}
\lim_{a\to0} a^{-n}\,P_n(x;a,a,1,1;q)=q^n\,h_n(x;q).
\end{equation*}

The quadratic transformations \eqref{39}, \eqref{40} immediately
imply quadratic transformations \cite[\S14.21]{14} connecting
discrete $q$-Hermite II polynomials \cite[\S14.29]{13}
and $q$-Laguerre polynomials \cite[\S14.21]{13}
because these two orthogonal polynomials can be expressed as
$i^{-n} h_n(ix;q^{-1})$ and
$\const p_n(-x;q^{-\al};q^{-1})$ in terms of
discrete $q$-Hermite I polynomials and little $q$-Laguerre polynomials,
respectively. Note that both families of orthogonal polynomials
have non-unique orthogonality measures, see for instance \cite{16}.
Quite probably these last quadratic transformations
are limit cases of rewritings of \eqref{37}, \eqref{38}
which can be interpreted as quadratic transformations for
pseudo big $q$-Jacobi polynomials \cite[Prop.~2.2]{15}.
\end{example}

Next we turn to the discrete part of the $q$-Askey scheme.
\begin{example}
On top there is a quadratic transformation between $q$-Racah
polynomials (see \cite[\S14.2]{13})
\begin{equation}
R_n(q^{-x}+\ga\de q^{x+1};\al,\be,\ga,\de\,|\,q)
:=\qhyp43{q^{-n},\al\be q^{n+1},q^{-x},\ga\de q^{x+1}}
{\al q,\be\de q,\ga q}{q,q}\quad(n=0,1,\ldots,N),
\label{44}
\end{equation}
where $\al q$ or $\be\de q$ or $\ga q$ is equal to $q^{-N}$.
It reads, with $N\in\{\thalf,1,\tfrac32,\ldots\}$,
\begin{multline}
R_{2n}(q^{-x-N-\half}-q^{x-N-\half};\al,\al,q^{-2N-2},-1\,|\,q)\\
=R_n(q^{-2x-2N-1}+q^{2x-2N-1};
\al^2,q^{-1},q^{-2N-2},q^{-2N-2}\,|\,q^2)\quad
(n=0,1,\ldots,[N+\thalf]).
\label{56}
\end{multline}
Indeed, as a function of $q^{-x-N-\half}-q^{x-N-\half}$
the polynomials on the \LHS\ of
\eqref{56} are orthogonal
on the points $q^{-x-N-\half}-q^{x-N-\half}$
($x=-N-\thalf,-N+\thalf,\ldots,N+\thalf$) with respect
to the weights
\[
(q^x+q^{-x})\,
\frac{(\al^2 q^2;q^2)_{x+N+\half}}{(q^2;q^2)_{x+N+\half}}\,
\frac{(\al^2 q^2;q^2)_{-x+N+\half}}{(q^2;q^2)_{-x+N+\half}}\,
\]
while the polynomials on the \RHS\ are orthogonal on the points
$(q^{-x-N-\half}-q^{x-N-\half})^2$
($x$ running over $-N-\thalf,-N+\thalf,\ldots,-\thalf$ or $0$)
with respect to the same weights.
These weights are positive if $-1<q\al<1$.

In terms of $q$-hypergeometric functions \eqref{56} can be written as
\begin{multline*}
\qhyp43{q^{-2n},\al^2 q^{2n+1},q^{-x-N-\half},-q^{x-N-\half}}
{q\al,-q\al,q^{-2N-1}}{q,q}\\
=\qhyp43{q^{-2n},\al^2 q^{2n+1},q^{-2x-2N-1},q^{2x-2N-1}}
{q^2\al^2,q^{-2N-1},q^{-2N}}{q^2,q^2},
\end{multline*}
which is the case $a=q^{-n}$, $b=q^n\al$, $c=q^{-x-N-\half}$,
$d=q^{x-N-\half}$ of \eqref{7}.

Similarly, from \eqref{9}, we have the quadratic transformation
\begin{multline}
R_{2n+1}(q^{-x-N-\half}-q^{x-N-\half};\al,\al,q^{-2N-2},-1\,|\,q)\\
=\frac{q^{-x-N-\half}-q^{x-N-\half}}{1-q^{-2N-1}}\,
R_n(q^{-2x-2N-1}+q^{2x-2N-1};\al^2,q,q^{-2N-2},q^{-2N-2}\,|\,q^2),
\label{57}
\end{multline}
where $N\in\{\thalf,1,\tfrac32,\ldots\}$ and $n=0,1,\ldots,[N]$.
Formula \eqref{57} can also be proved by orthogonality.

The special case $\al=0$ of \eqref{56} and \eqref{57}
gives quadratic transformations involving dual $q$-Krawtchouk
polynomials \cite[\S14.17]{13} and
dual $q$-Hahn polynomials \cite[\S14.7]{13}.
\end{example}

\begin{remark}
The quadratic transformations \eqref{37}, \eqref{38}
involving big and little $q$-Jacobi polynomials
can be obtained as limit cases of \eqref{56} and \eqref{57}.
For this we need the following special case of the
limit formula \cite[(2.2)]{18} from $q$-Racah polynomials to
big $q$-Jacobi polynomials:
\begin{equation*}
\lim_{N\to\iy}R_n(q^{-2N-1}x;a,a,q^{-2N-2},-1\,|\,q)
=\frac{P_n(x;a,a,1,1;q)}{P_n(-1;a,a,1,1;q)}\,.
\end{equation*}
We need also a limit formula from $q$-Racah polynomials to
little $q$-Jacobi polynomials, not yet observed in \cite{18}:
\begin{equation}
\lim_{N\to\iy}R_n(q^{-2N}x;a,b,q^{-N-1},\de q^{-N}\,|\,q)
=\frac{p_n(\de^{-1}x;b,a;q)}{p_n(1;b,a;q)}\,.
\end{equation}
This is obtained from the limit formula (straightforward from
\eqref{44})
\begin{equation*}
\lim_{N\to\iy}R_n(q^{-2N}x;a,b,q^{-N-1},\de q^{-N}\,|\,q)
=\qhyp31{q^{-n},abq^{n+1},\de x^{-1}}{qa}{q,b^{-1}\de^{-1}x}
\end{equation*}
combined with \cite[(III.8)]{3} and \eqref{36}.

Furthermore, the quadratic transformations \eqref{39}, \eqref{40}
can be obtained as limits of the cases $\al=0$ of \eqref{56}
and \eqref{57}.
\end{remark}

\begin{example}
Rather non-standard quadratic transformations for $q$-Racah polynomials
can be obtained by another specialization of
\eqref{7} and \eqref{9}:
\begin{align}
&\qhyp43{q^{-2n},q^{-2(N-n)-1},q^{-x},-\ga q^{x+1}}
{q^{-N},-q^{-N},\ga q}{q,q}=
\qhyp43{q^{-2n},q^{2n-2N-1},q^{-2x},\ga^2 q^{2x+2}}
{q^{-2N},\ga q,\ga q^2}{q^2,q^2},
\label{42}\\
&\qhyp43{q^{-2n-1},q^{-2(N-n)},q^{-x},-\ga q^{x+1}}{q^{-N},-q^{-N},\ga q}
{q,q}=\frac{q^{-x}-\ga q^{x+1}}{1-\ga q}\nonu\\
&\qquad\qquad\qquad\qquad\qquad\qquad\qquad\qquad\quad\times
\qhyp43{q^{-2n},q^{2n-2N+1},q^{-2x},\ga^2 q^{2x+2}}
{q^{-2N},\ga q^2,\ga q^3}{q^2,q^2}.
\label{43}
\end{align}
Here $N$ is a positive integer and $n=0,1,\ldots,N$.
For $2n\le N$ \eqref{42} is valid for all $x\in\CC$. However,
by the subtlety of passing to a lower parameter $q^{-N}$ in
\eqref{7} or \eqref{9}, formula \eqref{42} is only valid for $x=0,1,\ldots,N$ if $2n>N$.
Similarly, \eqref{43} is valid for all $x\in\CC$ if $2n+1\le N$,
but only valid for $x=0,1,\ldots,N$ if $2n+1>N$.

By substitution of \eqref{44} in \eqref{42} and \eqref{43}
we obtain quadratic transformations for $q$-Racah polynomials:
\begin{multline}
R_n(q^{-2x}+\ga^2 q^{2x+2};q^{-2N-2},q^{-1},\ga,\ga\,|\,q^2)\\
=\begin{cases}
R_{2n}(q^{-x}-\ga q^{x+1};q^{-N-1},q^{-N-1},\ga,-1\,|\,q)&(2n\le N),\\
R_{2N-2n+1}(q^{-x}-\ga q^{x+1};q^{-N-1},q^{-N-1},\ga,-1\,|\,q)&(2n>N),
\end{cases}
\label{45}
\end{multline}
\begin{multline}
\frac{q^{-x}-\ga q^{x+1}}{1-\ga q}\,R_n(q^{-2x}+\ga^2 q^{2x+2};q^{-2N-2},q,\ga,\ga\,|\,q^2)\\
=\begin{cases}
R_{2n+1}(q^{-x}-\ga q^{x+1};q^{-N-1},q^{-N-1},\ga,-1\,|\,q)&
(2n+1\le N),\\
R_{2N-2n}(q^{-x}-\ga q^{x+1};q^{-N-1},q^{-N-1},\ga,-1\,|\,q)&
(2n+1> N).
\end{cases}
\label{46}
\end{multline}
Both in \eqref{45} and \eqref{46} the identities corresponding to the
first case of the \RHS\ are valid for all complex
$y:=q^{-x}-\ga q^{x+1}$ (then
$q^{-2x}+\ga^2 q^{2x+2}=y^2+2\ga q$). But the identities corresponding to
the second case of the \RHS\ are only valid for $x=0,1,\ldots,N$.

By \cite[(14.2.2)]{13} the $q$-Racah polynomials on the \LHS\ of
\eqref{45} are orthogonal on the set of points
$q^{-2x}+\ga^2 q^{2x+2}$ ($x=0,1,\ldots,N$) with respect to the weights
\begin{equation}
w_x=q^{(2N+1)x}\,\frac{1+q^{2x+1}\ga}{1+q\ga}\,
\frac{(q^{-2N},q^2\ga^2;q^2)_x}{(q^2,q^{2N+4}\ga^2;q^2)_x}\,.
\label{47}
\end{equation}
These weights are positive if $q^{-N}<\ga<q^{-N-2}$.
Inspection of the positivity of the coefficient of $p_{n-1}(x)$ in
\cite[(14.2.4)]{13} for $n=1,\ldots,N$ gives the same constraint on
$\ga$.
Again by \cite[(14.2.2)]{13}, the $q$-Racah polynomials
on the \RHS\ of \eqref{45} are orthogonal on the set of points
$q^{-x}-\ga q^{x+1}$ ($x=0,1,\ldots,N$) with respect to the weights
$w_x$ given by \eqref{47}. This is compatible with \eqref{45},
but on the other hand \eqref{45} can be proved
from this equality of weights only if $2n\le N$.
Similar remarks can be made about \eqref{46}.

If we put for $n=0,1,\ldots,2N+1$
\[
p_n(y):=
\begin{cases}
R_{\half n}(y^2+2q\ga;q^{-2N-2},q^{-1},\ga,\ga\,|\,q^2)&\mbox{($n$ even),}\\
(1-\ga q)^{-1}y^{-1}\,
R_{\half(n-1)}(y^2+2q\ga;q^{-2N-2},q,\ga,\ga\,|\,q^2)&\mbox{($n$ odd)}
\end{cases}
\]
then $p_n(-y)=(-1)^n p_n(y)$ and the $p_n$ are orthogonal on
the set of points $\pm(q^{-x}-\ga q^{x+1})$ (for $x=0,1,\ldots,N$) with
respect to the weights $w_x$ given by \eqref{47}. For $n\le N$ the
explicit expressions for the $p_n$ as polynomials of general argument
are given by the first cases of the \RHS s of \eqref{45}, \eqref{46},
but the expressions for $n>N$ will be more complicated.
\end{example}

\subsection{Further examples of quadratic transformations in the
Askey scheme}
First we discuss limit cases for $q\uparrow1$
of the quadratic transformations
in the continuous part of the
Askey scheme.

\begin{example}
Between Wilson polynomials \cite[\S9.1]{13}
\begin{equation*}
\frac{W_n(x^2;a,b,c,d)}{W_n(-a^2;a,b,c,d)}=
\frac{W_n(x^2;a,b,c,d)}{(a+b)_n(a+c)_n(a+d)_n}
:=\hyp43{-n,a+b+c+d-1,a+ix,a-ix}{a+b,a+c,a+d}1
\end{equation*}
and continuous Hahn polynomials \cite[\S9.4]{13}
\begin{equation*}
\frac{p_n(x;a,b,\overline a,\overline b)}
{p_n(ia;a,b,\overline a,\overline b)}
=\frac{n!\,p_n(x;a,b,\overline a,\overline b)}
{i^n(a+\overline a)_n(a+\overline b)_n}
:=
\hyp32{-n,n+2\Re(a+b)-1,a+ix}{a+\overline a,a+\overline b}1
\end{equation*}
there are the quadratic transformations
\begin{align}
\frac{p_{2n}(x;a,b,\overline a,\overline b)}
{p_{2n}(ia;a,b,\overline a,\overline b)}&=
\frac{W_n(x^2;a,b,\thalf,0)}{W_n(-a^2;a,b,\thalf,0)}\,,
\label{49}\\
\frac{p_{2n+1}(x;a,b,\overline a,\overline b)}
{p_{2n+1}(ia;a,b,\overline a,\overline b)}&=
\frac{x W_n(x^2;a,b,\thalf,1)}{ia W_n(-a^2;a,b,\thalf,1)}\,,
\label{50}
\end{align}
where $a,b\in\RR$ or $b=\overline a$.
This follows by comparing the orthogonality relations
\cite[(9,1,2), (9.4.2)]{13} with each other.

In fact, \eqref{49} and \eqref{50} are limit cases of
the quadratic transformations \eqref{34}, \eqref{35}
for Askey-Wilson polynomials by the limits
\begin{equation*}
\lim_{q\uparrow 1}\frac{p_n(1-\thalf x(1-q)^2;q^a,q^b,q^c,q^d\,|\, q)}{(1-q)^{3n}}
=W_n(x;a,b,c,d)
\end{equation*}
and
\begin{multline*}
\lim_{q\uparrow1}
\frac{p_n\big(\cos\phi-x(1-q)\sin\phi;q^a e^{i\phi},q^b e^{i\phi},q^{\overline a} e^{-i\phi},
q^{\overline b} e^{-i\phi}\,|\, q\big)}
{(1-q)^{2n}}\\
=(-2\sin\phi)^n\,n!\,p_n(x;a,b,\overline a,\overline b)\qquad
(0<\phi<\pi).
\end{multline*}
There are corresponding limit cases of \eqref{7} and \eqref{9}:
\begin{align}
\hyp32{2a,2b+1,c}{a+b+1,c+d}1&=
\hyp43{a,b+\thalf,c,d}{a+b+1,\thalf(c+d),\thalf(c+d+1)}1,
\label{59}\\
\hyp32{2a,2b+1,c}{a+b+1,c+d}1&=\frac{d-c}{d+c}\,
\hyp43{a+\thalf,b+1,c,d}{a+b+1,\thalf(c+d)+1,\thalf(c+d+1)}1,
\label{58}
\end{align}
which are valid whenever both sides terminate.

Also note that \eqref{49} and \eqref{50} have the quadratic
transformations \eqref{10} as limit cases.
This follows by \cite[(9.4.15)]{13}
and the limit (extension of \cite[(9.1.18)]{13})
\begin{equation*}
\lim_{t\to\iy}
\frac{W_n(\thalf(1-x)t^2;a,\al+1-a,c+it,\be+1-c-it)}{t^{2n}n!}
=P_n^{(\al,\be)}(x).
\end{equation*}
\end{example}

\begin{example}
Between continuous dual Hahn polynomials \cite[\S9.3]{13}
\begin{equation*}
\frac{S_n(x^2;a,b,c)}{S_n(-a^2;a,b,c)}:=
\hyp32{-n,a+ix,a-ix}{a+b,a+c}1,\quad
S_n(-a^2;a,b,c)=(a+b)_n(a+c)_n,
\end{equation*}
and Meixner-Pollaczek polynomials \cite[\S9.7]{13}
\begin{equation*}
\frac{P_n^{(\la)}(x;\phi)}{P_n^{(\la)}(i\la;\phi)}:=
\hyp21{-n,\la+ix}{2\la}{1-e^{-2i\phi}},\quad
P_n^{(\la)}(i\la;\phi)=\frac{(2\la)_n}{n!}\,e^{in\phi},
\end{equation*}
there are the quadratic transformations
\begin{align}
\frac{P_{2n}^{(a)}(x;\thalf\pi)}{P_{2n}^{(a)}(ia;\thalf\pi)}
&=\frac{S_n(x^2;a,\thalf,0)}{S_n(-a^2;a,\thalf,0)}\,,
\label{51}\\
\frac{P_{2n+1}^{(a)}(x;\thalf\pi)}{P_{2n+1}^{(a)}(ia;\thalf\pi)}
&=\frac{x S_n(x^2;a,\thalf,1)}{ia S_n(-a^2;a,\thalf,1)}\,.
\label{52}
\end{align}
These are limit cases of \eqref{49} and \eqref{50} by the limits
\cite[(9.1.16), (9.4.14)]{13}
Furthermore, \eqref{51} and \eqref{52} have the quadratic
transformations \eqref{53} as limit cases by \cite[(9.7.15)]{13}
and the limit
\begin{equation}
\lim_{a\to\iy}\frac{S_n(ax;a,b,c)}{a^n n!}=L_n^{b+c-1}(x).
\label{54}
\end{equation}
For the proof of \eqref{54} compare the recurrence relations
\cite[(9.3.5), (9.12.4)]{13} with each other.
\end{example}

Next we turn to the discrete part of the Askey scheme.
\begin{example}
On top there is a quadratic transformation between Racah polynomials
\cite[\S9.2]{13}
\begin{multline*}
R_n\big(x(x+\ga+\de+1);\al,\be,\ga,\de\big):=
\hyp43{-n,n+\al+\be+1,-x,x+\ga+\de+1}{\al+1,\be+\de+1,\ga+1}1\\
(\mbox{$\al+1$ or $\be+\de+1$ or $\ga+1=-N$; } n=0,1,\ldots,N)
\end{multline*}
and Hahn polynomials \cite[\S9.5]{13}
\begin{equation*}
Q_n(x;\al,\be,N):=\hyp32{-n,n+\al+\be+1,-x}{\al+1,-N}1\quad
(n=0,1,\ldots,N).
\end{equation*}
It reads
\begin{multline}
Q_{2n}\big(x+N+\thalf;\al,\al,2N+1\big)
=R_n\big(x^2-(N+\thalf)^2;\al,-\thalf,-N-1,-N-1\big)\\
(N\in\{\thalf,1,\tfrac32,\ldots\},\; n=0,1,\ldots,[N+\thalf]).
\label{55}
\end{multline}
Indeed, as a function of $x$ the polynomials on the \LHS\ of
\eqref{55} are orthogonal
on the points $x=-N-\thalf,-N+\thalf,\ldots,N+\thalf$ with respect
to the weights
\[
\frac{(\al+1)_{N+\half+x}(\al+1)_{N+\half-x}}
{(N+\thalf+x)!\,(N+\thalf-x)!}\,,
\]
while the polynomials on the \RHS\ are orthogonal on the points $x^2$
($x$ running over $-N-\thalf,-N+\thalf,\ldots,-\thalf$ or $0$)
with respect to the same weights.

The quadratic transformation \eqref{55} is the case
$a=-n$, $b=n+\al$, $c=-x-N-\thalf$, $d=x-N-\thalf$
of formula \eqref{59}.
By specialization of \eqref{58}, also as a limit case for $q\uparrow1$
of \eqref{57}, we have the quadratic transformation
\begin{multline}
Q_{2n+1}(x+N+\thalf;\al,\al,2N+1)
=\frac{2N+1-2x}{2N+1}\,R_n(x^2-(N+\thalf)^2;\al+1,\thalf,-N-1,-N-1)\\
(N\in\{\thalf,1,\tfrac32,\ldots\},\; n=0,1,\ldots,[N]).
\label{60}
\end{multline}

The quadratic transformations \eqref{10} for Jacobi polynomials
can be obtained as limit cases of \eqref{55} and \eqref{60}.
\end{example}
\begin{example}
Quadratic transformations involving
Krawtchouk polynomials \cite[\S9.11]{13}
\begin{equation*}
K_n(x;p,N):=\hyp21{-n,-x}{-N}{p^{-1}}
\end{equation*}
and dual Hahn polynomials \cite[\S9.6]{13}
\begin{equation*}
R_n(x(x+\ga+\de+1);\ga,\de,N):=
\hyp32{-n,-x,x+\ga+\de+1}{\ga+1,-N}1
\end{equation*}
are given by
\begin{align}
K_{2m}(x+N;\thalf,2N)&=\frac{(\thalf)_m}{(-N+\thalf)_m}\,
R_m(x^2;-\thalf,-\thalf,N),
\label{71}\\
K_{2m+1}(x+N;\thalf,2N)&=-\,\frac{(\tfrac32)_m}{N\,(-N+\thalf)_m}\,
x\,R_m(x^2-1;\thalf,\thalf,N-1),
\label{72}\\
K_{2m}(x+N+1;\thalf,2N+1)&=\frac{(\tfrac12)_m}{(-N-\thalf)_m}\,
R_m(x(x+1);-\thalf,\thalf,N),
\label{73}\\
K_{2m+1}(x+N+1;\thalf,2N+1)&=\frac{(\tfrac32)_m}{(-N-\thalf)_{m+1}}\,
(x+\thalf)\,R_m(x(x+1);\thalf,-\thalf,N).
\label{74}
\end{align}
They can be proved by orthogonality, they are limit cases of
\eqref{55} and \eqref{60}, and they have the quadratic transformations
\eqref{53} involving Hermite and Laguerre polynomials as limit
cases.
\end{example}
\subsection{The ($q$-)Askey scheme of quadratic transformations}
Let us summarize the quadratic transformations for families in
the ($q$-)Askey scheme. In the $q$-case we have:
\begin{description}
\item[1a]
Askey-Wilson \eqref{34}, \eqref{35}
\item[1b]
$q$-Racah \eqref{56}, \eqref{57}
\item[2 \,]
big $q$-Jacobi to little $q$-Jacobi \eqref{37}, \eqref{38}
\item[3a]
Askey-Wilson \eqref{34}, \eqref{35} for $b=0$
\item[3b]
$q$-Racah \eqref{56}, \eqref{57} for $\al=0$
\item[4 \,]
discrete $q$-Hermite to Wall \eqref{39}, \eqref{40}
\item[5 \,]
Askey-Wilson \eqref{34}, \eqref{35} for $a=b=0$
\end{description}
The transformations 1a and 1b in $q$-hypergeometric form are related
by analytic continuation, similarly for 3a and 3b.
The limit arrows between the various cases are as follows.
\[
\begin{matrix}
{\rm 1a}&&&&{\rm 1b}\\
\downarrow&\searrow&&\swarrow&\downarrow\\
{\rm 3a}&&2&&{\rm 3b}\\
\downarrow&\searrow&\downarrow&\swarrow&\\
5&&4&&
\end{matrix}
\]

In the case $q=1$ we have:
\begin{description}
\item[1a]
continuous Hahn to Wilson \eqref{49}, \eqref{50}
\item[1b]
Hahn to Racah \eqref{55}, \eqref{60}
\item[2 \,]
Jacobi \eqref{10}
\item[3a]
Meixner-Pollaczek to continuous dual Hahn \eqref{51}, \eqref{52}
\item[3b]
Krawtchouk to dual Hahn \eqref{71}--\eqref{74}
\item[4 \,]
Hermite to Laguerre \eqref{53}
\end{description}
The transformations 1a and 1b in hypergeometric form are related
by analytic continuation, similarly for 3a and 3b.
The limit arrows between the various cases are as above, except
that the case 5 is missing. There are limits for $q\uparrow1$
from the $q$-cases to the corresponding $q=1$ cases. The $q$-case 5
also has a limit to the $q=1$ case 4.
\section{The two-variable case}
\subsection{General polynomials}
For an analogue of \eqref{1} in two variables we generalize the proof
of Theorem 10.1 in Sprinkhuizen~\cite{4}.
We will work with monomials $x^{m-l}y^l$ ($m,l\in\ZZ$, $m\ge l\ge0$)
with a dominance partial ordering
\[
(m,l)\le(n,k)\quad{\rm iff}\quad m\le n\quad{\rm and}\quad m+l\le n+k.
\]
Let $w(x,y)$ be a (nonnegative) weight function on a domain $\Om\subset\RR^2$ such that
\[
\int_\Om |x|^{m-l} |y|^l w(x,y)\,dx\,dy<\iy\quad\mbox{for all $m,l$.}
\]
Let $p_{n,k}(x,y)$ be polynomials of the form
\begin{equation}
p_{n,k}(x,y)=\sum_{(m,l)\le(n,k)}c_{m,l}x^{m-l}y^l,\quad c_{n,k}\ne0,
\label{11}
\end{equation}
such that
\begin{equation}
\int_\Om p_{n,k}(x,y)\,x^{m-l}y^l\,w(x,y)\,dx\,dy=0\quad{\rm if}\quad
(m,l)<(n,k).
\label{12}
\end{equation}
We call the polynomials $p_{n,k}(x,y)$ {\em dominance orthogonal}
polynomials.
For convenience we assume that they are monic, i.e., $c_{n,k}=1$ in
\eqref{11}.
Thus $p_{n,k}(x,y)$ and $p_{m,l}(x,y)$ with $(n,k)\ne(m,l)$ are
orthogonal on $\Om$ with weight function $w(x,y)$ if
$(n,k)$ and $(m,l)$ are related in the partial ordering $\le$,
but the orthogonality will usually fail if $(n,k)$ and $(m,l)$ are
not related in this partial ordering, except for very special $\Om$ and
$w(x,y)$, as will occur for cases related to root systems.

Now suppose that $\Om$ is invariant under $(x,y)\to(-x,y)$, and also
$w(x,y)=w(-x,y)$.
Then, by \eqref{12},
$p_{n,k}(-x,y)=(-1)^{n-k} p_{n,k}(x,y)$, and in \eqref{11} $c_{m,l}=0$
if $n-k$ and $m-l$ do not have the same parity.

Put
\begin{equation}
\Om':=\{(y,x^2)\mid(x,y)\in\Om\}\quad{\rm and}\quad
v(x,y):=w(y^\half,x)\quad((x,y\in\Om').
\label{14}
\end{equation}
\begin{proposition}
\label{17}
Let $q_{n,k}(x,y)$ and $r_{n,k}(x,y)$ be dominance orthogonal polynomials
on $\Om'$ with respect to weight functions $y^{-\half}v(x,y)$ and
$y^\half v(x,y)$, respectively. Then
\begin{equation}
q_{n,k}(y,x^2)=p_{n+k,n-k}(x,y),\quad x\,r_{n,k}(y,x^2)=p_{n+k+1,n-k}(x,y).
\label{13}
\end{equation}
\end{proposition}
\Proof
We have
\[
p_{n+k,n-k}(x,y)=\sum_{(i,j)\le(n+k,n-k)}c_{i,j}x^{i-j}y^j,
\]
where only terms with $i-j$ even occur. So we can substitute $i-j=2l$
and $i+j=2m$, $c_{i,j}=c'_{m,l}$. Then $(i,j)\le(n+k,n-k)$ iff $(m,l)\le(n,k)$. Hence
\[
p_{n+k,n-k}(x,y)=\sum_{(m,l)\le(n,k)}c'_{m,l}y^{m-l} x^{2l},
\]
while from \eqref{12} we have
\begin{equation*}
\int_{\Om'} p_{n+k,n-k}(y^\half,x)\,x^{m-l}y^l\,v(x,y)y^{-\half}\,dx\,dy=0\quad{\rm if}\quad
(m,l)<(n,k).
\end{equation*}
This settles \eqref{13} for $q_{n,k}$. A similar proof can be given
for $r_{n,k}$\,.\qed
\bPP
In particular, let $\Om$ be the region
\begin{equation}
\Om:=\{(x,y)\in\RR^2\mid 1-x+y,1+x+y,x^2-4y>0\}.
\label{15}
\end{equation}
Then \eqref{14} and \eqref{15} yield that
\begin{equation}
\Om'=\{(x,y)\in\RR^2\mid y,y-4x,(1+x)^2-y>0\}
=\{(\thalf x,1+x+y)\mid (x,y)\in\Om\}.
\label{16}
\end{equation}
So if $\Om$ is given by \eqref{15} then, by \eqref{16}, an affine
transformation respecting the dominance partial order of monomials
maps $\Om'$ onto $\Om$. Thus we can formulate a variant of Proposition
\ref{17} which again generalizes the proof
of Theorem 10.1 in Sprinkhuizen~\cite{4}:
\begin{proposition}
\label{20}
Let $\Om$ be given by \eqref{15}.
Let the $p_{n,k}(x,y)$ be monic dominance orthogonal polynomials on $\Om$
with respect to a weight function $w(x,y)=w(-x,y)$.
Define $v(x,y)$ on $\Om$ by
\begin{equation*}
v(2y,x^2-2y-1)=w(x,y).
\end{equation*}
Let $q_{n,k}(x,y)$ and $r_{n,k}(x,y)$ be dominance orthogonal polynomials
on $\Om$ with respect to weight functions $(1+x+y)^{-\half}v(x,y)$ and
$(1+x+y)^\half v(x,y)$, respectively. Then
\begin{align}
2^{-n+k}q_{n,k}(2y,x^2-2y-1)&=p_{n+k,n-k}(x,y),
\label{18}\\
2^{-n+k}x r_{n,k}(2y,x^2-2y-1)&=p_{n+k+1,n-k}(x,y).
\label{19}
\end{align}
\end{proposition}

If the $p_{n,k}(x,y)$ in the above Proposition are not monic
but satisfy $p_{n,k}(2,1)\ne0$ (which probably is implied by the
dominance orthogonality)
then we can replace \eqref{18}, \eqref{19} by
\begin{align}
\frac{q_{n,k}(2y,x^2-2y-1)}{q_{n,k}(2,1)}&=
\frac{p_{n+k,n-k}(x,y)}{p_{n+k,n-k}(2,1)}\,,
\label{28}\\
\frac{x r_{n,k}(2y,x^2-2y-1)}{2r_{n,k}(2,1)}&=
\frac{p_{n+k+1,n-k}(x,y)}{p_{n+k+1,n-k}(2,1)}\,.
\label{29}
\end{align}

In further variants of these results, to be discussed below,
we will formulate results in a normalization as in \eqref{28},
\eqref{29}. If the assumption corresponding to $p_{n,k}(2,1)\ne0$
would fail then formulations in terms of monic polynomials would
still be true.
\subsection{Symmetric polynomials}
In Proposition \ref{20} replace $x,y$ by $\xi,\eta$, and next put $\xi=x+y$,
$\eta=xy$. Then we can rephrase this proposition in terms of symmetric
polynomials in $x,y$. For this purpose make the following observations.
\begin{itemize}
\item
The map $(x,y)\to(\xi,\eta)$ is a diffeomorphism from
\begin{equation}
\La:=\{(x,y)\mid -1<y<x<1\}
\label{26}
\end{equation}
onto $\Om$ given by \eqref{15}.
Furthermore $d\xi\,d\eta=(x-y)\,dx\,dy$.
\item
Let $n>k$. Then, for certain $a_i,b_i$ with $a_0=b_0=1$ we have
\begin{align*}
(x+y)^{n-k}(xy)^k&=
\sum_{i=0}^{[\half(n-k)]}a_i(x^{n-i}y^{k+i}+x^{k+i}y^{n-i}),\\
x^ny^k+x^ky^n&=
\sum_{i=0}^{[\half(n-k)]}b_i(x+y)^{n-k-2i}(xy)^{k+i}.
\end{align*}
\item
If
\begin{equation}
p(\xi,\eta)=\sum_{(m,l)\le(n,k)} a_{m,l} \xi^{m-l}\eta^l
\label{21}
\end{equation}
for certain $a_{m,l}$ with $a_{n,k}\ne0$ then
\begin{equation}
p(x+y,xy)=\sum_{(m,l)\le(n,k)} b_{m,l}(x^my^l+x^ly^m)
\label{22}
\end{equation}
for certain $b_{m,l}$ with $b_{n,k}\ne0$. Conversely, any symmetric
polynomial given by the \RHS\ of \eqref{22} can be written as $p(x+y,xy)$
for some polynomial $p(\xi,\eta)$ of the form \eqref{21}.
\end{itemize}

Now let $W(x,y)$ be a weight function on $\La$ and let $P_{n,k}(x,y)$ be
symmetric polynomials of the form of the \RHS\ of
\eqref{22} with $b_{n,k}\ne0$ such that
\begin{equation*}
\int_\La P_{n,k}(x,y) (x^my^l+x^ly^m)W(x,y)(x-y)\,dx\,dy=0\quad
{\rm if}\quad (m,l)<(n,k).
\end{equation*}
We call the polynomials $P_{n,k}(x,y)$ dominance orthogonal
symmetric polynomials. Observe that the polynomials $p_{n,k}(\xi,\eta)$
are dominance orthogonal on $\Om$ with weight function $w(\xi,\eta)$
iff the polynomials $P_{n,k}(x,y):=p_{n,k}(x+y,xy)$ are dominance orthogonal
on $\La$ with orthogonality measure $W(x,y)(x-y)\,dx\,dy$, where
$W(x,y):=w(x+y,xy)$.

Now we can rephrase Proposition \ref{20} as follows.
\begin{proposition}
\label{25}
Let the $P_{n,k}(x,y)$ be dominance orthogonal symmetric polynomials
on $\La$
with respect to a measure $W(x,y)(x-y)\,dx\,dy$, where
$W(x,y)=W(-y,-x)$.
Define a weight function $V$ on $\La$ by
\begin{equation}
V(xy+(1-x^2)^\half(1-y^2)^\half,xy-(1-x^2)^\half(1-y^2)^\half)=W(x,y).
\label{27}
\end{equation}
Let $Q_{n,k}(x,y)$ and $R_{n,k}(x,y)$ be dominance orthogonal symmetric
polynomials
on $\La$ with respect to measures
$(1+x)^{-\half}(1+y)^{-\half}V(x,y)(x-y)\,dx\,dy$ and
$(1+x)^{\half}(1+y)^{\half}V(x,y)(x-y)\,dx\,dy$, respectively. Then
\begin{align}
\frac{Q_{n,k}(xy+(1-x^2)^\half(1-y^2)^\half,xy-(1-x^2)^\half(1-y^2)^\half)}
{Q_{n,k}(1,1)}&=
\frac{P_{n+k,n-k}(x,y)}{P_{n+k,n-k}(1,1)}\,,
\label{23}\\
\frac{(x+y) R_{n,k}(xy+(1-x^2)^\half(1-y^2)^\half,
xy-(1-x^2)^\half(1-y^2)^\half)}{2R_{n,k}(1,1)}&=
\frac{P_{n+k+1,n-k}(x,y)}{P_{n+k+1,n-k}(1,1)}\,.
\label{24}
\end{align}
if $P_{n,k}(1,1)\ne0$ for all $n,k$,
or the same identities without denominators for monic polynomials.
\end{proposition}

On passing to trigonometric coordinates \eqref{23} and \eqref{24}
can be rewritten as
\begin{align}
\frac{Q_{n,k}(\cos(\tha_1-\tha_2),\cos(\tha_1+\tha_2))}
{Q_{n,k}(1,1)}&=
\frac{P_{n+k,n-k}(\cos\tha_1,\cos\tha_2)}{P_{n+k,n-k}(1,1)}\,,
\\
\frac{(\cos\tha_1+\cos\tha_2)
R_{n,k}(\cos(\tha_1-\tha_2),\cos(\tha_1+\tha_2))}{2R_{n,k}(1,1)}&=
\frac{P_{n+k+1,n-k}(\cos\tha_1,\cos\tha_2)}{P_{n+k+1,n-k}(1,1)}\,.
\end{align}
\begin{example}
\label{62}
In the notation of Proposition \ref{20} let
\begin{equation*}
w(\xi,\eta)=w_{\al,\be,\ga}(\xi,\eta):=
(1-\xi+\eta)^\al(1+\xi+\eta)^\be (\xi^2-4\eta)^\ga\quad(\al,\be,\ga>-1,\;
\al+\ga,\be+\ga>-\tfrac32),
\end{equation*}
and put $p_{n,k}(\xi,\eta)=p_{n,k}^{\al,\be,\ga}(\xi,\eta)$
for the corresponding
dominance orthogonal polynomials on the region $\Om$ defined by \eqref{15}.
These polynomials, nowadays known as Jacobi polynomials for root system
$BC_2$, were first studied in \cite{5} and subsequently in \cite{4}.
It follows from \cite[(3.14)]{5} and \cite[Theorem 8.1]{4}
that these polynomials, even if they are defined as dominance orthogonal
polynomials, still satisfy full orthogonality, and that they are
nonzero at $(2,1)$ by the explicit value \cite[(7.3)]{4}.

Since
\begin{equation*}
w_{\ga,0,\al}(2\eta,\xi^2-2\eta-1)=4^\al w_{\al,\al,\ga}(\xi,\eta),
\end{equation*}
we have
\begin{equation}
\frac{p_{n,k}^{\ga,-\half,\al}(2\eta,\xi^2-2\eta-1)}
{p_{n,k}^{\ga,-\half,\al}(2,1)}=\frac{p_{n+k,n-k}^{\al,\al,\ga}(\xi,\eta)}
{p_{n+k,n-k}^{\al,\al,\ga}(2,1)}\,,\quad
\frac{\xi\,p_{n,k}^{\ga,\half,\al}(2\eta,\xi^2-2\eta-1)}
{2p_{n,k}^{\ga,\half,\al}(2,1)}=
\frac{p_{n+k+1,n-k}^{\al,\al,\ga}(\xi,\eta)}
{p_{n+k+1,n-k}^{\al,\al,\ga}(2,1)}\,.
\label{64}
\end{equation}
These quadratic transformations were first given by Sprinkhuizen
\cite[Theorem 10.1]{4}. They can be conceptually explained
by the fact that $B_2$ and $C_2$, while special cases of $BC_2$, are
isomorphic root systems.

Equivalently, in the notation of Proposition \ref{25}, let
\begin{equation*}
W(x,y)=W_{\al,\be,\ga}(x,y):=
(1-x)^\al(1-y)^\al(1+x)^\be(1+y)^\be(x-y)^{2\ga}
\end{equation*}
and put $P_{n,k}(x,y)=P_{n,k}^{\al,\be,\ga}(x,y)$ for the corresponding
dominance orthogonal symmetric polynomials on the region $\La$. Then
\begin{align*}
\frac{P_{n,k}^{\ga,-\half,\al}(\cos(\tha_1-\tha_2),\cos(\tha_1+\tha_2))}
{P_{n,k}^{\ga,-\half,\al}(1,1)}&=
\frac{P_{n+k,n-k}^{\al,\al,\ga}(\cos\tha_1,\cos\tha_2)}
{P_{n+k,n-k}^{\al,\al,\ga}(1,1)}\,,
\\
\frac{(\cos\tha_1+\cos\tha_2)
P_{n,k}^{\ga,\half,\al}(\cos(\tha_1-\tha_2),\cos(\tha_1+\tha_2))}
{2P_{n,k}^{\ga,\half,\al}(1,1)}&=
\frac{P_{n+k+1,n-k}^{\al,\al,\ga}(\cos\tha_1,\cos\tha_2)}
{P_{n+k+1,n-k}^{\al,\al,\ga}(1,1)}\,.
\end{align*}
\end{example}
\subsection{Symmetric Laurent polynomials}
Let $S_2$ be the symmetric group in $2$ letters and
$W_2:=S_2\ltimes(\ZZ_2)^2$ (the Weyl group of $BC_2$).
These groups naturally act on $\ZZ^2$.
For $\la=(\la_1,\la_2)\in\ZZ^2$ and $x=(x_1,x_2)\in\CC^2$
put $x^\la:=x_1^{\la_1}x_2^{\la_2}$.
Put
\begin{equation*}
m_\la(x):=\sum_{\mu\in S_2\la} x^\mu,\qquad
\wt m_\la(z):=\sum_{\mu\in W_2\la} z^\mu\qquad
(\la_1\ge\la_2\ge0).
\end{equation*}
For certain $a_\mu,b_\mu$ with $a_\la,b_\la\ne0$ we have
\begin{align*}
m_\la\big(\thalf(z_1+z_1^{-1}),\thalf(z_2+z_2^{-1})\big)
&=\sum_{\mu\le\la} a_\mu \wt m_\mu(z_1,z_2),\\
\wt m_\la(z_1,z_2)&=\sum_{\mu\le\la} b_\mu
m_\mu\big(\thalf(z_1+z_1^{-1}),\thalf(z_2+z_2^{-1})\big).
\end{align*}

Let $W(x,y)$ be a weight function on the region $\La$ given by \eqref{26}.
Let
\begin{equation}
\Ga:=\big\{(z_1,z_2)\in\CC^2\;\big|\; |z_1|=|z_2|=1,\;
0<\arg z_1<\arg z_2<\pi\big\}.
\label{31}
\end{equation}
Then
$(z_1,z_2)\mapsto\big(\thalf(z_1+z_1^{-1}),\thalf(z_2+z_2^{-1})\big)$
is a diffeomorphism from $\Ga$ onto $\La$.
On $\Ga$ define a weight function $\De(z_1,z_2)$ such that
\begin{equation*}
W\big(\thalf(z_1+z_1^{-1}),\thalf(z_2+z_2^{-1})\big)
=\frac{8\De(z_1,z_2)}{(z_1-z_1^{-1})(z_2-z_2^{-1})
(z_1+z_1^{-1}-z_2-z_2^{-1})}\,.
\end{equation*}
Then
\begin{equation*}
\int_\De f(x,y)\,W(x,y)\,(x-y)\,dx\,dy
=\int_\Ga f\big(\thalf(z_1+z_1^{-1}),\thalf(z_2+z_2^{-1})\big)\,
\De(z_1,z_2)\,\frac{dz_1}{z_1}\,\frac{dz_2}{z_2}\,.
\end{equation*}
Hence, if the $P_{n,k}(x,y)$ are dominance orthogonal symmetric polynomials
on $\La$ with orthogonality measure $W(x,y)(x-y)\,dx\,dy$ and if
\begin{equation*}
p_{n,k}(z_1,z_2):=P_{n,k}\big(\thalf(z_1+z_1^{-1}),\thalf(z_2+z_2^{-1})\big)
\end{equation*}
then the $p_{n,k}(z_1,z_2)$ are dominance orthogonal $W_2$-invariant
Laurent polynomials on $\Ga$ with weight function $\De(z_1,z_2)$,
i.e., we have for certain $c_{m,l}$ with $c_{n,k}\ne0$ that
\begin{equation}
p_{n,k}(z_1,z_2)=\sum_{(m,l)\le(n,k)}c_{m,l}\,\wt m_{m,l}(z_1,z_2)
\label{32}
\end{equation}
such that
\begin{equation*}
\int_\De p_{n,k}(z_1,z_2)\,\wt m_{m,l}(z_1,z_2)\,\De(z_1,z_2)\,
\frac{dz_1}{z_1}\,\frac{dz_2}{z_2}=0\quad{\rm if}\quad(m,l)<(n,k).
\end{equation*}
Call the polynomials $p_{n,k}(z_1,z_2)$ monic if $c_{n,k}=1$ in
\eqref{32}.

Now we can rephrase Proposition \ref{25} as follows.
\begin{proposition}
Let the $p_{n,k}(z_1,z_2)$ be dominance orthogonal $W_2$-invariant
polynomials on $\Ga$
with respect to a weight function $\De(z_1,z_2)$, where
$\De(z_1,z_2)=\De(-z_2^{-1},-z_1^{-1})$.
Define a weight function $\wt\De$ on $\Ga$ by
\begin{equation*}
\wt\De(z_1z_2,z_1z_2^{-1})=\De(z_1,z_2).
\end{equation*}
Let $q_{n,k}(x,y)$ and $r_{n,k}(x,y)$ be dominance orthogonal $W_2$-invariant
polynomials on $\Ga$ with respect to weight functions
$\wt\De(z_1,z_2)$ and
$(1+z_1)(1+z_1^{-1})(1+z_2)(1+z_2^{-1})\wt\De(z_1,z_2)$, respectively. Then
\begin{align}
\frac{q_{n,k}(z_1z_2,z_1z_2^{-1})}
{q_{n,k}(1,1)}&=
\frac{p_{n+k,n-k}(z_1,z_2)}{p_{n,k}(1,1)}\,,\\
\frac{(z_1+z_1^{-1}+z_2+z_2^{-1}) r_{n,k}(z_1z_2,z_1z_2^{-1})}
{4r_{n,k}(1,1)}&=
\frac{p_{n+k+1,n-k}(z_1,z_2)}{p_{n+k+1,n-k}(1,1)}\,.
\end{align}
if $p_{n,k}(1,1)\ne0$ for all $n,k$, or the same identities without
denominators for monic polynomials.
\end{proposition}

\begin{example}
In the notation of Proposition \ref{20} let
\begin{equation*}
\De(z_1,z_2)=\De(z_1,z_2;q,t;a,b,c,d)=
\De_+(z_1,z_2)\De_+(z_1^{-1},z_2^{-1}),
\end{equation*}
where
\begin{equation*}
\De_+(z_1,z_2):=\frac{(z_1^2;q)_\iy}
{(az_1,bz_1,cz_1,dz_1;q)_\iy}\,
\frac{(z_2^2;q)_\iy}
{(az_2,bz_2,cz_2,dz_2;q)_\iy}
\,\frac{(z_1z_2,z_1z_2^{-1};q)_\iy}
{(tz_1z_2,tz_1z_2^{-1};q)_\iy},
\end{equation*}
and put $p_{n,k}(z_1,z_2)=p_{n,k}(z_1,z_2;q,t;a,b,c,d)$
for the corresponding dominance orthogonal $W_2$-invariant monic
Laurent polynomials on the region $\Ga$ defined by \eqref{31}.
These polynomials are the two-variable case of the $n$-variable Koornwinder polynomials \cite{6}, \cite{7}, which are
associated with root system $BC_n$. These polynomials are fully
orthogonal \cite{6}. Now observe that
\begin{equation*}
\De(z_1,z_2;q,t;a,-a,q^\half,-q^\half)
=\De(z_1z_2,z_1z_2^{-1};q^2,a^2;t,qt,-1,-q).
\end{equation*}
Hence, by Proposition \ref{20}, we have for $n+k$ even that
\begin{align}
P_{n,k}(z_1,z_2;q,t;a,-a,q^\half,-q^\half)&=
P_{\half(n+k),\half(n-k)}(z_1z_2,z_1z_2^{-1};q^2,a^2;t,qt,-1,-q),
\label{69}\\
P_{n+1,k}(z_1,z_2;q,t;a,-a,q^\half,-q^\half)&=(z_1+z_2+z_1^{-1}+z_2^{-1})
\nonu\\
&\quad\times
P_{\half(n+k),\half(n-k)}(z_1z_2,z_1z_2^{-1};q^2,a^2;t,qt,-q,-q^2).
\label{70}
\end{align}
\end{example}
\subsection{Failure of quadratic transformations in the $n$-variable case
if $n>2$}
There are no straightforward analogues in $n>2$ variables of
Propositions \ref{20} and \ref{25}.
Indeed, symmetric polynomials in $x_1,\ldots,x_n$ invariant under
$x_i\to-x_i$ ($i=1,\ldots,n$) correspond to polynomials in
$e_1,\ldots,e_n$ (the elementary symmetric polynomials in
$x_1,\ldots,x_n$) which are invariant under $e_{2i-1}\to-e_{2i-1}$
($i=1,\ldots,[\thalf(n+1)]$). If $n>2$ then this last involutive
linear transformation has more than one eigenvalue unequal to 1.
Therefore, by Stanley \cite[Theorem 4.1]{8} (a theorem going back
to Shephard \& Todd \cite{9}), there do not exist $n$ algebraically
independent invariants for this involution if $n>2$.
\section{Discussion of results and further perspective}
\subsection{New results suggested by extrapolation from very few data}
In all examples of quadratic transformations within multi-parameter
families of special orthogonal polynomials in one or two variables
we start with a subfamily depending on less than the full number
of parameters, and then
there is an even degree case and an odd degree case giving rise
to two systems of orthogonal polynomials for which one of the
parameters takes two special values, say $-\half$ in the even case
and $\half$ in the odd case. Thus formulas and other results already
known for the system with which we started give results for these
parameter values $\pm\half$ which can be tentatively extrapolated
for more general values of the parameter.

\begin{example}
Consider the quadratic tranformations 
\eqref{10} for Jacobi polynomials.
They map from the Gegenbauer case of parameters $(\al,\al)$
to the Jacobi cases $(\al,\pm\half)$.
The Gegenbauer case is easier than the general Jacobi case,
so \eqref{10} may be helpful as a start to derive from known
results in the Gegenbauer case yet unknown results in the
general Jacobi case.
For instance, for a system of orthogonal polynomials $\{p_n\}$
it is remarkable to have a lowering formula of the form
\begin{equation*}
\frac d{dx}\Big(\psi(x)^n\,p_n(\phi(x))\Big)
=\la_n \psi(x)^{n-1}\,p_{n-1}(\phi(x)).
\end{equation*}
For Gegenbauer polynomials such a formula does exist
(see \cite[(3.3)]{19}):
\begin{equation}
\frac d{dx}\left((1+x^2)^{\half n}\
P_n^{(\al,\al)}\left(\frac x{\sqrt{1+x^2}}\right)\right)
=(n+\al)\,(1+x^2)^{\half(n-1)}\,
P_{n-1}^{(\al,\al)}\left(\frac x{\sqrt{1+x^2}}\right),
\label{61}
\end{equation}
but probably not for general Jacobi polynomials.
But let us see what we get for $(\al,\pm\half)$ by quadratic
transformation of \eqref{61}. First we have to iterate \eqref{61} once.
Then apply \eqref{10}. We obtain
\begin{align*}
\frac{d^2}{dx^2}\left((1+x^2)^n\,
P_n^{(\al,-\half)}\left(\frac{x^2-1}{x^2+1}\right)\right)
&=4(n+\al)(n-\thalf)\,(1+x^2)^{n-1}\,
P_{n-1}^{(\al,-\half)}\left(\frac{x^2-1}{x^2+1}\right),\\
\left(\frac{d^2}{dx^2}+\frac2x\,\frac d{dx}\right)\left((1+x^2)^n\,
P_n^{(\al,\half)}\left(\frac{x^2-1}{x^2+1}\right)\right)
&=4(n+\al)(n+\thalf)\,(1+x^2)^{n-1}\,
P_{n-1}^{(\al,\half)}\left(\frac{x^2-1}{x^2+1}\right).
\end{align*}
Then the straightforward extrrapolation
\begin{equation*}
\left(\frac{d^2}{dx^2}+\frac{2\be+1}x\,\frac d{dx}\right)\left((1+x^2)^n\,
P_n^{(\al,\be)}\left(\frac{x^2-1}{x^2+1}\right)\right)
=4(n+\al)(n+\be)\,(1+x^2)^{n-1}\,
P_{n-1}^{(\al,\be)}\left(\frac{x^2-1}{x^2+1}\right)
\end{equation*}
can indeed be proved, see \cite[(4.4)]{19}.
\end{example}
\begin{example}
\label{68}
In this example we again have a result obtaiend by quadratic
transformation, now valid on a two-dimensional subdomain
of a three-dimensional parameter space,
but still can make a meaningful guess
how to extrapolate.

The $BC_2$-type Jacobi polynomials $p_{n,k}^{\al,\be,\ga}(\xi,\eta)$
given in Example \ref{62} have an explicit expansion
\cite[(6.11)]{20} (with different notation following
\cite[(3.3)]{20})
in terms of polynomials
\begin{equation}
(1-\xi+\eta)^{\half(m+l)}\,
P_{m-l}^{(\ga,\ga)}\left(\frac{1-\thalf\xi}{(1-\xi+\eta)^\half}\right)
\quad(0\le l\le m\le n,\;l\le k).
\label{63}
\end{equation}
The polynomials \eqref{63} can be recognized as Jack polynomials
in two variables and the mentioned expansion was seen in
\cite[Section 11.2]{7} as a limit case of Okounkov's binomial
formula for Koornwinder polynomials in two variables.
Now by \eqref{64} and parity we have a quadratic transformation
\begin{equation}
p_{n+k,n-k}^{\al,\al,\ga}(\xi,\eta)=\const
p_{n,k}^{-\half,\ga,\al}(-2\eta,\xi^2-2\eta-1).
\label{65}
\end{equation}
We can explicitly expand the \RHS\ of \eqref{65} in terms of
polynomials
\begin{equation}
\xi^{m+l} P_{m-l}^{(\al,\al)}\left(\frac{1+\eta}\xi\right),
\label{66}
\end{equation}
and thus this is also an explicit expansion for the \LHS\
of \eqref{65}. The other quadratic transformation in \eqref{64}
gives a similar result for $p_{n+k+1,n-k}^{\al,\al,\ga}(\xi,\eta)$
(with $x^{m+l}$ in \eqref{66} replaced by $x^{m+l+1}$).
This suggests, and is indeed confirmed in \cite[Section 7]{20},
that $p_{n,k}^{\al,\be,\ga}(\xi,\eta)$ has a nice expansion in
terms of the polynomials
\begin{equation}
\xi^m P_l^{(\al,\be)}\left(\frac{1+\eta}\xi\right)\quad
(0\le m-l\le n-k,\;m+l\le n+k),
\label{67}
\end{equation}
which can be considered as a two-parameter extension of the
Jack polynomials in two variables.

In fact, by \cite[Theorem 7.7]{20},
the polynomials \eqref{63} and \eqref{67} are limit cases of
$p_{n,k}^{\al,\be,\ga}(\xi,\eta)$ for $\beta\to\iy$ and
$\ga\to\iy$, respectively.
\end{example}
\subsection{Further perspective}
In the one-variable part of this paper we gave a quite extensive
treatment of quadratic transformations between families in the
Askey and $q$-Askey scheme. Similar treatments should be given
in the two-variable case. On the one hand we have orthogonal
polynomials in two variables which are products of two
polynomials from the ($q$-)Askey scheme and an elementary function,
of which the orthogonal polynomials on the triangle involving
products of two Jacobi polynomials are a well-known example.
Quadratic transformations for such polynomials can be derived 
by suitable substitutions of quadratic transformations for polynomials
in one variable. On the other hand there are the orthogonal polynomials
associated with root system $BC_2$. By work of various authors
a large part of the ($q$-)Askey scheme has now been realized for $BC_2$.
It can be expected that corresponding schemes of quadratic
transformations can also be given in the $BC_2$ case.

Finally it would be interesting to do further
explicit work for Koornwinder polynomials in two variables
analogous to the $q=1$ case treated in \eqref{20} and extending
\cite[Section 11.1]{7}. Analogous to Example \ref{68} for $q=1$,
the quadratic transformations \eqref{69}, \eqref{70}
may be helpful for making a start in such work.
\quad\\
\begin{footnotesize}
\begin{quote}
{ T. H. Koornwinder, Korteweg-de Vries Institute, University of
 Amsterdam,\\
 P.O.\ Box 94248, 1090 GE Amsterdam, The Netherlands;

\vspace{\smallskipamount}
email: }{\tt T.H.Koornwinder@uva.nl}
\end{quote}
\end{footnotesize}


\begin{thebibliography}{99}
%
%until 20
%
\bibitem{1}
R. Askey and J. A. Wilson,
{\em Some basic hypergeometric orthogonal polynomials that generalize Jacobi polynomials},
Mem. Amer. Math. Soc. 54 (1985), no. 319.
%
\bibitem{16}
C. Berg,
{\em On some indeterminate moment problems for measures on a geometric
progression},
J. Comput. Appl. Math. 99 (1998), 67--75.
%
\bibitem{10}
T. S. Chihara,
{\em An introduction to orthogonal polynomials},
Gordon and Breach, 1978;
reprinted, Dover, 2011.
%
\bibitem{3}
G.~Gasper and M.~Rahman,
{\em Basic hypergeometric series}, 2nd edn.,
Cambridge University Press, 2004.
%
\bibitem{15}
W. Groenevelt and E. Koelink,
{\em The indeterminate moment problem for the $q$-Meixner polynomials},
J. Approx. Theory 163 (2011), 836--863.
%
\bibitem{13}
R. Koekoek, P. A. Lesky and R. F. Swarttouw,
{\em Hypergeometric orthogonal polynomials and their $q$-analogues},
Springer-Verlag, 2010.
%
\bibitem{5}
T. H. Koornwinder,
{\em Orthogonal polynomials in two variables which are eigenfunctions
of two algebraically independent partial differential operators, I, II},
Nederl. Akad. Wetensch. Proc. Ser. A 77 = Indag. Math. 36 (1974),
48--58, 59--66.
%
\bibitem{6}
T. H. Koornwinder,
{\em Askey-Wilson polynomials for root systems of type $BC$},
in: {\em Hypergeometric functions on domains of positivity,
Jack polynomials, and applications}, Contemp. Math. 138,
American Mathematical Society, 1992, pp.~189--204.
%
\bibitem{11}
T. H. Koornwinder,
{\em Askey-Wilson polynomials as zonal spherical functions on the $SU(2)$ quantum group}, SIAM J. Math. Anal. 24 (1993), 795--813.
%
\bibitem{19}
T. H. Koornwinder,
{\em Lowering and raising operators for some special orthogonal polynomials},
in: {\em Jack, Hall-Littlewood and Macdonald polynomials},
Contemp. Math 417, American Mathematical Society, 2006, pp. 227--238;
{\tt arXiv:math/0505378v1 [math.CA]}. 
%
\bibitem{18}
T. H. Koornwinder,
{\em On the limit from q-Racah polynomials to big q-Jacobi polynomials}, SIGMA 7 (2011), 040, 8 pp.;
{\tt arXiv:1011.5585v4 [math.CA]}.
%
\bibitem{12}
T. H. Koornwinder,
{\em $q$-Special functions, a tutorial},
{\tt arXiv:math/9403216v2 [math.CA]}, 2013.
%
\bibitem{14}
T.~H. Koornwinder,
{\em Additions to the formula lists in ``Hypergeometric orthogonal
polynomials and their q-analogues'' by Koekoek, Lesky and Swarttouw},  {\tt arXiv:1401.0815v2 [math.CA]}, 2015.
%
\bibitem{7}
T. H. Koornwinder,
{\em Okounkov's BC-type interpolation Macdonald polynomials and their
$q=1$ limit},
S\'em. Lothar. Combin. B72a (2015), 27 pp.;
{\tt arXiv:1408.5993v5 [math.CA]}.
%
\bibitem{20}
T. H. Koornwinder and I. G. Sprinkhuizen,
{\em Generalized power series expansions for a class of orthogonal
polynomials in two variables},
SIAM J. Math. Anal. 9 (1978), 457--483.
%
\bibitem{9}
G. C. Shephard and J. A. Todd,
{\em Finite unitary reflection groups},
Canadian J. Math. 6 (1954), 274--304.
%
\bibitem{2}
V. N. Singh,
{\em The basic analogues of identities of the Cayley-Orr type},
J. London Math. Soc. 34 (1959), 15--22.
%
\bibitem{4}
I. G. Sprinkhuizen-Kuyper,
{\em Orthogonal polynomials in two variables. A further analysis of the
polynomials orthogonal over a region bounded by two lines and a parabola},
SIAM J. Math. Anal. 7 (1976), 501--518.
%
\bibitem{8}
R. Stanley,
{\em Invariants of finite groups and their applications to combinatorics},
Bull. Amer. Math. Soc. (N.S.) 1 (1979), 475--511.
%
\bibitem{17}
G. Szeg{\H{o}},
{\em Orthogonal polynomials},
Colloquium Publications 23,
American Mathematical Society, Fourth Edition, 1975.
%
\end{thebibliography}
\end{document}